\documentclass[a4paper,12pt,doublespace]{article}
\usepackage{epsfig,graphicx}
\usepackage{psfrag}
\usepackage{amssymb,bbm}
\usepackage{amsmath}
\usepackage{amsfonts}
\usepackage{setspace}
\usepackage[toc,page]{appendix}
\usepackage{enumerate}
\usepackage{tikz}

\usepackage{fancyhdr}
\usepackage{pstricks}
\usepackage{amsfonts}
\usepackage{pst-infixplot}
\usepackage{pst-math}
\usepackage{sidecap}

\usetikzlibrary{patterns,shapes,snakes,backgrounds,calc}
\tikzset{axis line style/.style={thin, gray, -stealth}}
\usepackage{float}
\usepackage{tikz-cd}
\usepackage{blindtext}
\usepackage{upgreek}
\usepackage{fancyhdr}
\usepackage{subcaption}

\usepackage{pstricks,pst-math,pst-plot}
\usepackage[dotinlabels]{titletoc}
\usepackage{psfrag}
\usepackage{pstricks-add}

\usepackage[colorlinks]{hyperref}
\newcounter{theorem}

\newtheorem{definition}{Definition}

\title{A Weighted-Graph Curvature Calculator and Whether the Discrete Curvature Senses the Smooth One\date{}}
\author{G\"{o}k\c{c}e \c{C}AKMAK \footnote{Eskisehir Technical University, Science Faculty, Department of Mathematics, 26470, Eski\c{s}ehir, Turkey.  e-mail: gokcecakmak@eskisehir.edu.tr}
\thanks{Corresponding Author.}
\and Ali DEN\.{I}Z \footnote{Eskisehir Technical University, Science Faculty, Department of Mathematics, 26470, Eski\c{s}ehir, Turkey. e-mail: adeniz@eskisehir.edu.tr } \and  \c{S}ahin KO\c{C}AK \footnote{Anadolu University (emeritus), Science Faculty, Department of Mathematics, 26470, Eski\c{s}ehir, Turkey. e-mail: skocak@anadolu.edu.tr} \and Murat L\.IMONCU \footnote{Eskisehir Technical University, Science Faculty, Department of Mathematics, 26470, Eski\c{s}ehir, Turkey.  e-mail: mlimoncu@eskisehir.edu.tr}}
\begin{document}
\maketitle
\thispagestyle{empty}

\begin{abstract}
We investigate whether there is a relationship between the discrete Bakry-\'{E}mery curvature of a graph and the smooth curvature of an ambient surface into which the graph is embedded geodesically. As we used weighted graphs as test objects, we developed a program for the calculation of the discrete curvature and with the help of this calculator, we observed some indications of such a relationship.
\end{abstract}


\section{Introduction}
The notion of graph curvature has attracted considerable interest in recent years (\cite{CushingCalc2022}, \cite{CushingEig2022}, \cite{Fathi2022}, \cite{Klartag2016}, \cite{Yau2011}, \cite{Viola2021}). There are several approaches to this discrete curvature notion and one of them, the so-called Bakry-\'{E}mery graph curvature, has direct roots in smooth differential geometry (\cite{Bauer2017}). We wondered whether the Bakry-\'{E}mery curvature of a graph somehow reflects the smooth curvature of an ambient surface in which the graph is embedded in such a way that the edges are realized as geodesic segments on the surface. To investigate and test this idea, it would be natural and appropriate to consider weighted graphs. As computations of discrete graph curvatures by hand are notoriously difficult and the existing main graph curvature calculator \cite{CushingCalc2022} does not include weighted graphs, we have developed a calculator for weighted graphs (see Appendix).

To test the above hope of reflection of the smooth surface curvature in the discrete curvature of geodesically embedded graphs, we chose as a test-object the following weighted "umbrella graph" $G=(V, E, \omega)$ with the vertex set $V=\{v_0, v_1, \dots , v_n\}$, the edge set $E=\{[v_0v_i], [v_iv_{i+1}]\; | \; i=1,\dots,n,\; mod\, n\}$ and the weight function $\omega:E\to \mathbb{R}^+$, $\omega([v_0v_i])=1,$ $\omega([v_iv_{i+1}])=\rho$ with some fixed $\rho$. For a certain value of $\rho$ (depending on $n$), the umbrella graph can be geodesically embedded into the Euclidean plane and below and above this value of $\rho$ (within a certain range) it can be embedded into a sphere or a hyperbolic plane with a certain curvature depending on $\rho$. We were eager to know how the discrete Bakry-\'{E}mery curvature changes in dependance of $\rho$.

We now fix the definitions and notations \cite{CushingEig2022}:
Given a finite, simple, connected, weighted graph $G=(V, E, \omega)$, the weighted Laplacian $\Delta= \Delta_\omega$ acting on $f:V\to \mathbb{R}$ is defined by
\[ \Delta_\omega f(x)=\dfrac{1}{\omega(x)}\displaystyle\sum_{v\sim x}\omega_{xv}(f(v)-f(x)),  \]
where $v\sim x$ denotes that $v$ is a neighbour of $x$ in $G$, $\omega_{xv}=\omega([xv])$ and $\omega(x)=\displaystyle\sum_{v\sim x}\omega_{xv}$.

The Laplacian gives rise to the symmetric bilinear forms
\begin{align*}
  2\Gamma_{\omega}(f,g) &:=\Delta_\omega(f,g)-f\Delta_\omega g-g\Delta_\omega f , \\
   2\Gamma_{\omega,2}(f,g)  & :=\Delta_\omega(\Gamma_{\omega}(f,g))-\Gamma_{\omega}(f, \Delta_\omega g)-\Gamma_{\omega}(g, \Delta_\omega f).
\end{align*}

$\Gamma_{\omega}(f,g)$ can be expressed in more explicit terms as
\[ \Gamma_{\omega}(f,g)(x)=\dfrac{1}{2\omega(x)}\displaystyle\sum_{v\sim x}\omega_{xv}(f(v)-f(x))(g(v)-g(x)). \]
We will use the abbreviations $\Gamma_{\omega}(f):=\Gamma_{\omega}(f,f)$ and $\Gamma_{\omega,2}(f):=\Gamma_{\omega, 2}(f,f)$.

\begin{definition}
The Bakry-\'{E}mery curvature for the dimension $\infty$ at a vertex $x\in V$ of a weighted graph $G=(V, E, \omega)$ is the maximum value $K\in \mathbb{R} \cup \{-\infty\}$ such that for any real function $f:V\to \mathbb{R}$, $\Gamma_{\omega,2}(f)(x)\geq K\Gamma_{\omega}(f)(x)$.
\end{definition}

We remark that the Bakry-\'{E}mery curvature $K_G(x)$ of the graph $G$ at the vertex $x$ can also be expressed as
\[
K_G(x)=\inf_{f} \dfrac{\Gamma_{\omega,2}(f)(x)}{\Gamma_{\omega}(f)(x)}, \qquad \qquad (\Gamma_{\omega}(f)(x)\neq 0).
\]
This can be shown along the lines of \cite{Viola2021} and it results from the fact that for $\Gamma_{\omega}(f)(x)=0$, the quantity $\Gamma_{\omega,2}(f)(x)$ is nonnegative. For computations, one can assume without loss of generality $f(x)=0$.

In the next section, we explain how to compute the Bakry-\'{E}mery curvature at a vertex and derive our program from that description. In the last section, we apply this device to the above mentioned umbrella graphs and observe promising indications of a relationship between the discrete and smooth Bakry-\'{E}mery curvatures as documented by the Figures \ref{grafik1}, \ref{grafik2}, \ref{grafik3}, \ref{grafik4} and Table \ref{table1}.

\section{The Computation of the Bakry-\'{E}mery Curvature}

We first compute the terms $\Gamma_{\omega}(f)(x)$ and $\Gamma_{\omega,2}(f)(x)$ for a function $f:V\to \mathbb{R}$ with $f(x)=0$. In the following, we will drop the subscript $\omega$ as we will stay always in the weighted setting. Obviously,
\begin{eqnarray*}
(\Delta f)(x)&=&\frac{1}{\omega(x)} \sum_{v\sim x} \omega_{xv} f(v),\\
\Gamma(f)(x)&=&\frac{1}{2\omega(x)} \sum_{v \sim x} \omega_{xv} f^2(v). \\
\end{eqnarray*}

We now elaborate $\Gamma_2(f)(x)$. First, we compute $\Delta\Gamma(f)(x)$:
\begin{align*}
  \Delta\Gamma(f)(x) & =\dfrac{1}{\omega(x)}\sum_{v \sim x} \omega_{xv} \left[ \Gamma(f)(v) -\Gamma(f)(x) \right]   \\
  & = \dfrac{1}{\omega(x)}\sum_{v \sim x} \omega_{xv} \Gamma(f)(v)-\dfrac{1}{\omega(x)}\Gamma(f)(x)\sum_{v \sim x}\omega_{xv}   \\
   & =\dfrac{1}{\omega(x)}\sum_{v \sim x} \omega_{xv}\dfrac{1}{2\omega(v)}\sum_{u \sim v} \omega_{uv}(f(u)-f(v))^2-\dfrac{1}{2\omega(x)}\sum_{v \sim x}\omega_{xv}f^2(v)  \\
   & =\dfrac{1}{\omega(x)}\sum_{v \sim x}\sum_{u \sim v}\dfrac{\omega_{xv}\omega_{uv}}{2\omega(v)}(f(u)-f(v))^2-\dfrac{1}{2\omega(x)}\sum_{v \sim x}\omega_{xv}f^2(v).
\end{align*}

Now, we compute $\Gamma(f,\Delta f)(x)$:
\begin{align*}
  2\Gamma(f,\Delta f)(x) & =2\dfrac{1}{2\omega(x)}\sum_{v \sim x} \omega_{xv} f(v)[\Delta f(v)-\Delta f(x)] \\
   & =\dfrac{1}{\omega(x)}\sum_{v \sim x} \omega_{xv} f(v)\left[ \dfrac{1}{\omega(v)}\sum_{u \sim v}\omega_{uv}(f(u)-f(v))-\dfrac{1}{\omega(x)}\sum_{z \sim x} \omega_{zx} f(z) \right] \\
   & =\dfrac{1}{\omega(x)}\sum_{v \sim x}\left[\omega_{xv} f(v) \dfrac{1}{\omega(v)}\sum_{u\sim v}\omega_{uv} f(u)-\omega_{xv} f(v) \dfrac{1}{\omega(v)}\sum_{u\sim v}\omega_{uv} f(v) \right]\\
   &\hspace{.4cm}-\dfrac{1}{\omega^2(x)}\sum_{z \sim x} \omega_{zx} f(z)\sum_{v \sim x} \omega_{xv} f(v) \\
   & =\dfrac{1}{\omega(x)}\sum_{v \sim x}\sum_{u \sim v}\dfrac{\omega_{xv}\omega_{uv}}{\omega(v)}f(u)f(v)-\dfrac{1}{\omega(x)}\sum_{v \sim x}\omega_{xv}f^2(v)\\
   &\hspace{.4cm}-\dfrac{1}{\omega^2(x)}\left( \sum_{v \sim x} \omega_{xv} f(v)\right)^2.
\end{align*}

We now insert these expressions into $\Gamma_2(f)(x)$:
\begin{align*}
  2\Gamma_2(f)(x) & =\Delta\Gamma(f)(x) - 2\Gamma(f,\Delta f)(x)\\
   & =\dfrac{1}{\omega(x)}\sum_{v \sim x}\sum_{u \sim v}\dfrac{\omega_{xv}\omega_{uv}}{2\omega(v)}(f(u)-f(v))^2-\dfrac{1}{2\omega(x)}\sum_{v \sim x}\omega_{xv}f^2(v) \\
   &\hspace{.4cm} -\dfrac{1}{\omega(x)}\sum_{v \sim x}\sum_{u \sim v}\dfrac{\omega_{xv}\omega_{uv}}{\omega(v)}f(u)f(v)+\dfrac{1}{\omega(x)}\sum_{v \sim x}\omega_{xv}f^2(v)\\
   &\hspace{.4cm}+\dfrac{1}{\omega^2(x)}\left( \sum_{v \sim x} \omega_{xv} f(v)\right)^2 \\
   & =\dfrac{1}{\omega(x)}\sum_{v \sim x}\sum_{u \sim v}\dfrac{\omega_{xv}\omega_{uv}}{2\omega(v)}(f(u)-f(v))^2+\dfrac{1}{2\omega(x)}\sum_{v \sim x}\omega_{xv}f^2(v)\\
   &\hspace{.4cm} -\dfrac{1}{\omega(x)}\sum_{v \sim x}\sum_{u \sim v}\dfrac{\omega_{xv}\omega_{uv}}{\omega(v)}f(u)f(v)+\dfrac{1}{\omega^2(x)}\left( \sum_{v \sim x} \omega_{xv} f(v)\right)^2\\
    & =\dfrac{1}{\omega(x)}\sum_{v \sim x}\sum_{u \sim v}\dfrac{\omega_{xv}\omega_{uv}}{\omega(v)}\left[\dfrac{1}{2}(f(u)-f(v))^2-f(u)f(v)\right]\\
    &\hspace{.4cm} +\dfrac{1}{2\omega(x)}\sum_{v \sim x}\omega_{xv}f^2(v) +\dfrac{1}{\omega^2(x)}\left( \sum_{v \sim x} \omega_{xv} f(v)\right)^2.
\end{align*}

\begin{figure}[H]
  \centering
  \includegraphics[scale=0.75]{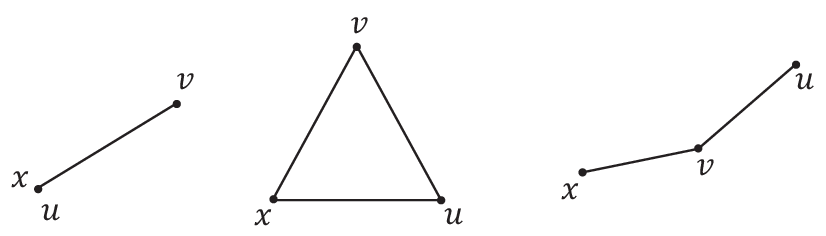}
   \caption{Here the first figure represents the case $x=u$, the second figure the triangular position and the third figure the linear position.} \label{sekil1}
\end{figure}

We decompose the pair of relationships $v\sim x$ and $u\sim v$ into three constellations (see Figure \ref{sekil1}):
\begin{enumerate}
  \item $u=x.$
  \item $u\sim x$ i.e. $u$ is a neighbour of $x$ and this case is often called triangular position. We denote the set of ordered triples $(x, v, u)$ in this position by the notation $T(x, v, u)$:
      \[ T(x, v, u)=\{(x, v, u) \; | \; v\sim x, u\sim v, u\sim x\}.\]
  \item Otherwise (combinatorial distance of $u$ to $x$ is two). For this case, which we call the ``linear" position, we denote the set of ordered triples $(x, v, u)$ in this position by the notation $L(x, v, u)$:
      \[ L(x, v, u)=\{(x, v, u) \; | \; v\sim x, u\sim v, u \mbox{ is not a neighbour of }x\}.\]
\end{enumerate}

We evaluate the first summand $\displaystyle S=\dfrac{1}{\omega(x)}\sum_{v \sim x}\sum_{u \sim v}\dfrac{\omega_{xv}\omega_{uv}}{\omega(v)}\left[\dfrac{1}{2}(f(u)-f(v))^2-f(u)f(v)\right]$ separately for the three cases stated above.
\begin{itemize}
  \item[1.] $S$ for $u=x$: $\displaystyle S_{u=x}=\dfrac{1}{\omega(x)}\sum_{v\sim x}\dfrac{\omega^2_{xv}}{2\omega(v)}f^2(v)$.
  \item[2.] $S$ for $T(x, v, u)$: $ \displaystyle S_{T(x, v, u)}=\dfrac{1}{\omega(x)}\sum_{T(x, v, u)}\dfrac{\omega_{xv}\omega_{uv}}{\omega(v)}\left[\dfrac{1}{2}(f(u)-f(v))^2-f(u)f(v)\right]$.
  \item[3.] $S$ for $L(x, v, u)$: $\displaystyle S_{L(x, v, u)}=\dfrac{1}{\omega(x)}\sum_{L(x, v, u)}\dfrac{\omega_{xv}\omega_{uv}}{\omega(v)}\left[\dfrac{1}{2}(f(u)-f(v))^2-f(u)f(v)\right]$.
\end{itemize}

We have a disjoint sum $S=S_{u=x}+S_{T(x, v, u)}+S_{L(x, v, u)}$. We now get
\[
   2\Gamma_2(f)(x)  = S_{L(x, v, u)}+S_{T(x, v, u)}+\dfrac{1}{2\omega(x)}\sum_{v\sim x}\left(\omega_{xv}+\dfrac{\omega^2_{xv}}{\omega(v)}\right)f^2(v)+\dfrac{1}{\omega^2(x)}\left( \sum_{v \sim x} \omega_{xv} f(v)\right)^2.
  \]

As we want to compute the infimum of $\dfrac{\Gamma_2(f)}{\Gamma(f)}$, we can get rid of the terms coming from vertices $u$ which have a combinatorial distance 2 to $x$.

Denote the subset of $L(x, v, u)$ with such a fixed vertex $u$ by $L^u(x, v, u)$. Then the part $S_{L^u(x, v, u)}$ of $S_{L(x, v, u)}$ containing the terms associated with $u$ is
\[ S_{L^u(x, v, u)}= \dfrac{1}{\omega (x)}\sum_{L^u(x, v, u)}\dfrac{\omega_{xv}\omega_{uv}}{\omega(v)}\left[\dfrac{1}{2}(f(u)-f(v))^2-f(u)f(v)\right].\]
Note that this sum runs over $v$ (by fixed $u$) which are intermediary vertices between $x$ and $u$.

The quadratic expression in $f(u)$ is minimized by
\begin{equation*}
 \displaystyle f(u)=\dfrac{\displaystyle 2\sum_{L^u(x, v, u)}\dfrac{\omega_{xv}\omega_{uv}}{\omega(v)}f(v)}{\displaystyle \sum_{L^u(x, v, u)}\dfrac{\omega_{xv}\omega_{uv}}{\omega(v)}}.
\end{equation*}

Inserting these values into $\Gamma_2(f)$ we get a quadratic expression in terms of $f(v)$ where the vertices $v$ are the neighbours of $x$. $\dfrac{\Gamma_2(f)}{\Gamma(f)}$ can then be infimized by standard linear algebra as the least eigenvalue of the matrix $2\omega(x)(w_{xv}^{-1})(\Gamma_2(f))^q(x)$ where $(w_{xv}^{-1})$ is a diagonal matrix and $(\Gamma_2(f))^q$ is the matrix associated with the quadratic form $\Gamma_2(f)(x)$ (cf. \cite{Viola2021}). We could in fact write this quadratic form explicitly, but the formidable formulas might not be necessary to be written down at this point. In the appendix, we give the source code of a calculator of the Bakry-\'{E}mery curvature of a weighted graph. In the next section we use this program to compute the curvature of the umbrella graphs defined in the introduction.

\section{The Curvature of the Umbrella Graphs}

\begin{figure}[H]
  \centering
  \includegraphics[scale=0.6]{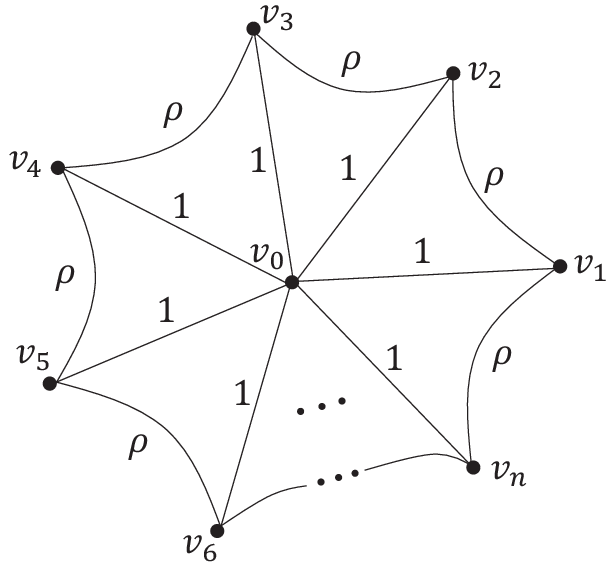}
   \caption{The umbrella graph $G_{n, \rho}$} \label{sekil2}
\end{figure}

We consider the following weighted graph $G_{n, \rho}=(V,E,\omega)$ with the vertex set $V=\{v_0, v_1, \dots , v_n\}$, the edge set $E=\{[v_0v_i], [v_iv_{i+1}]\; | \; i=1,\dots,n,\; mod\, n\}$ and the weight function $\omega :E\to \mathbb{R}^+$, $\omega([v_0v_i])=1, \omega([v_iv_{i+1}])=\rho$ with some fixed $\rho$ depicted in Figure \ref{sekil2}.

For $\rho= \sqrt{2-2\cos\frac{2\pi}{n}}$, this graph can be geodesically embedded into the Euclidean plane, for $0<\rho<\sqrt{2-2\cos\frac{2\pi}{n}}$ into the 2-sphere with appropriate radius, and for $\sqrt{2-2\cos\frac{2\pi}{n}}<\rho<2$ into the hyperbolic plane with an appropriate curvature. These radii and curvatures can be determined by the following spherical and hyperbolical cosine formulas: In spherical case
\[
\cos \dfrac{\rho}{R}=\cos^2\dfrac{1}{R}+\sin^2\dfrac{1}{R}\cos\dfrac{2\pi}{n}
\]
where $1/R^2$ is the curvature of the sphere, and in the hyperbolic case
\[
\cosh \dfrac{\rho}{K}=\cosh^2\dfrac{1}{K}-\sinh^2\dfrac{1}{K}\cos\dfrac{2\pi}{n}
\]
where $-1/K^2$ is the curvature of the hyperbolic plane.

Specifically, the values of $\rho$ for $R=1$ and $K=1$ will be used in the sequel. For a fixed $n$, we denote the value of $\rho$ for which the $n$-umbrella graph $G_{n,\rho}$ can be geodesically embedded into the sphere with radius 1 by $\rho_n^+$ and we denote the value of $\rho$ for which $G_{n,\rho}$ can be embedded into the hyperbolic plane with curvature -1 by $\rho_n^-$ (sometimes we denote the Euclidean embedding value $\sqrt{2-2\cos\frac{2\pi}{n}}$ of $\rho$ for $G_{n,\rho}$ by $\rho_n^0$).

\textbf{The case for n=3:} Let $n=3$ and $f:V=\{v_0,v_1,v_2, v_3\}\to \mathbb{R}$ be a function with $f(v_0)=0$. We denote the values $f(v_i)$ by $\alpha_i$ for $i=1,2,3$. Then
\[\Gamma(f)(v_0)=\dfrac{1}{6}\sum_{i=1}^{3}\alpha_i^2,\]
and
\[
\Gamma_2(f)(v_0)=\left(\dfrac{\rho}{6(1+2\rho)}+\dfrac{2}{9}\right)\sum_{i=1}^{3}\alpha_i^2+\dfrac{1}{9}\sum_{i<j}\alpha_i\alpha_j-\dfrac{2\rho}{3(1+2\rho)}
\sum_{\substack{i=1\\ mod 4}}^3\alpha_i\alpha_{i+1}.
\]

The infimum of $\dfrac{\Gamma_2(f)(v_0)}{\Gamma(f)(v_0)}$ is the least eigenvalue of the matrix
\[
\setlength\arraycolsep{8pt}
\begin{pmatrix}
 \dfrac{4}{3}+\dfrac{\rho}{1+2\rho} & \dfrac{1}{3}-\dfrac{2\rho}{1+2\rho} & \dfrac{1}{3}-\dfrac{2\rho}{1+2\rho} \\[2em]
  \dfrac{1}{3}-\dfrac{2\rho}{1+2\rho}  & \dfrac{4}{3}+\dfrac{\rho}{1+2\rho} & \dfrac{1}{3}-\dfrac{2\rho}{1+2\rho} \\[2em]
  \dfrac{1}{3}-\dfrac{2\rho}{1+2\rho}  & \dfrac{1}{3}-\dfrac{2\rho}{1+2\rho}  & \dfrac{4}{3}+\dfrac{\rho}{1+2\rho}
\end{pmatrix}.\]

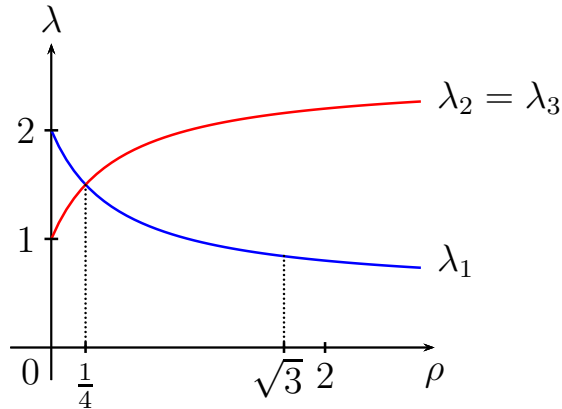
\begin{figure}[H]
\centering
\scalebox{1.2}{\begin{pspicture}*(-1,-1)(5.8,4.5)
\psset{xunit=1.5cm,yunit=1.2cm}
\SpecialCoor%
\psaxes[labels=none,ticks=none,linewidth=0.5pt]{->}(0,0)(-0.3,-0.3)(2.8,2.8)
\psPlot[linecolor=blue]{0}{2.7}{(2+x)/(1+2*x)}
\psPlot[linecolor=red]{0}{2.7}{(1+5*x)/(1+2*x)}
\uput[r](2.7,0.8){$\lambda_1$}
\uput[r](2.7,2.3){$\lambda_2=\lambda_3$}

\psline[linestyle=dotted,dotsep=0.7pt](0.25,0)(0.25,1.5)
\psline[linestyle=dotted,dotsep=0.7pt](1.7,0)(1.7,0.85)
\uput[l](0,1){1}  \psline(-0.05,1)(0.05,1)
\uput[l](0,2){2}  \psline(-0.05,2)(0.05,2)
\uput[d](0.25,0){$\frac{1}{4}$} \psline(0.25,-0.05)(0.25,0.05)
\uput[dl](0,0){0} \uput[d](2,0){2} \psline(2,-0.05)(2,0.05)
\uput[d](2.8,0){$\rho$} \uput[d](0,3.3){$\lambda$}
\uput[d](1.65,0.05){$\sqrt{3}$} \psline(1.7,-0.05)(1.7,0.05)
\end{pspicture}}
\caption{The graph of the eigenvalues for $G_{3,\rho}$.}{\label{grafik1}}
\end{figure}

The eigenvalues are $\lambda_1=\dfrac{2+\rho}{1+2\rho}$ and $\lambda_2=\lambda_3=\dfrac{1+5\rho}{1+2\rho}$, with the least one being $\lambda_1$ on the interval $[1/4, \infty)$ (see Figure \ref{grafik1}) and $\lambda_2=\lambda_3$ on the interval $(0, 1/4]$.

\begin{figure}[H]
  \centering
  \includegraphics[scale=0.6]{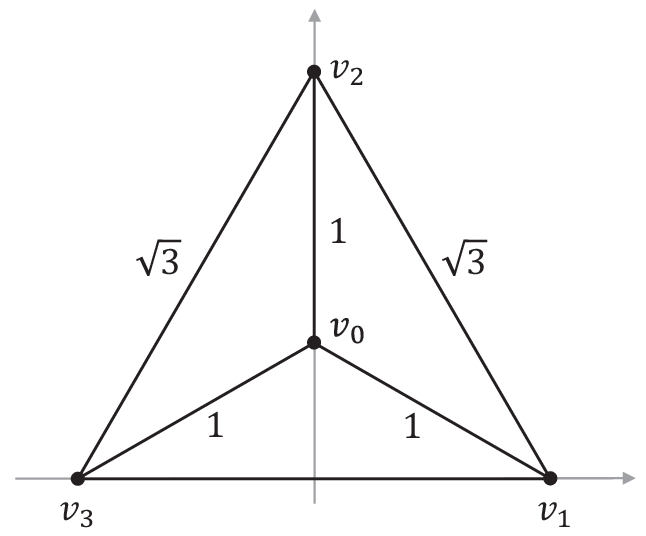}
   \caption{The umbrella graph $G_{3, \rho}$ with $\rho=\sqrt{3}$.} \label{sekil7}
\end{figure}

For $\rho=\sqrt{3}$, the umbrella graph $G_{3, \rho}$ is geodesically embeddable into the Euclidean plane (see Figure \ref{sekil7}) with Bakry-\'{E}mery curvature $K_G(v_0)=\dfrac{2+\sqrt{3}}{1+2\sqrt{3}}$.

\begin{figure}[H]
  \centering
  \includegraphics[scale=0.6]{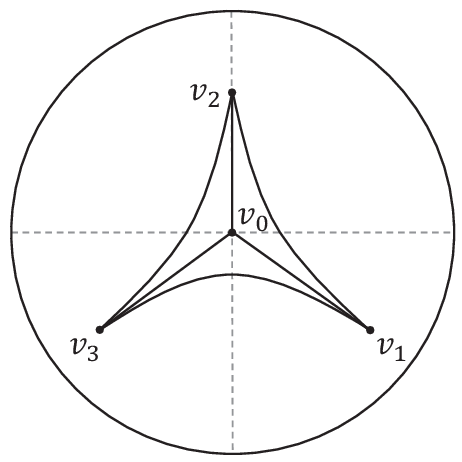}
   \caption{The umbrella graph $G_{3, \rho}$ with $\sqrt{3}<\rho<2$.} \label{sekil8}
\end{figure}

For $\sqrt{3}<\rho<2$, $G_{3, \rho}$ is geodesically embeddable into the hyperbolic plane with (smooth) curvature $-\dfrac{1}{K^2}$ (which is the solution of the equation $\cosh \dfrac{\rho}{K}=\cosh^2\dfrac{1}{K}+\dfrac{1}{2}\sinh^2\dfrac{1}{K}$ by hyperbolic trigonometry) (see Figure \ref{sekil8}) and has the discrete curvature $K_G(v_0)=\dfrac{2+\rho}{1+2\rho}$.

\begin{figure}[H]
  \centering
  \includegraphics[scale=0.6]{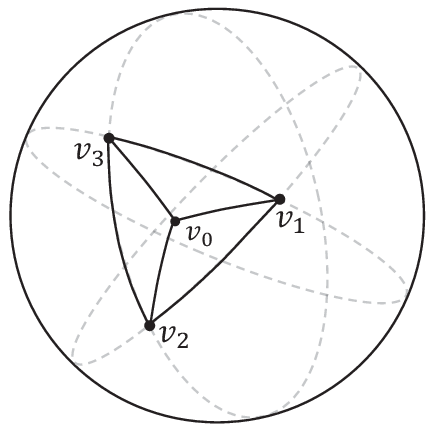}
   \caption{The umbrella graph $G_{3, \rho}$ with $0<\rho<\sqrt{3}$.} \label{sekil9}
\end{figure}

For $0<\rho<\sqrt{3}$, the graph $G_{3, \rho}$ is geodesically embeddable into the 2-sphere with radius $R$ (which is the solution of the equation $\cos \dfrac{\rho}{R}=\cos^2\dfrac{1}{R}-\dfrac{1}{2}\sin^2\dfrac{1}{R}$ by spherical trigonometry) and has the discrete curvature $\dfrac{2+\rho}{1+2\rho}$ for $1/4<\rho<\sqrt3$ and $\dfrac{1+5\rho}{1+2\rho}$ for $0<\rho<1/4$ (see Figure \ref{sekil9}).

\textbf{The case for n=4:} Let $n=4$ and $f:V=\{v_0,v_1,\dots, v_4\}\to \mathbb{R}$ be a function with $f(v_0)=0$. We denote the values $f(v_i)$ by $\alpha_i$ for $i=1,\dots,4$. Then
\[\Gamma(f)(v_0)=\dfrac{1}{8}\sum_{i=1}^{4}\alpha_i^2,\]
and
\[
\Gamma_2(f)(v_0)=\left(\dfrac{\rho}{8(1+2\rho)}+\dfrac{5}{32}\right)\sum_{i=1}^{4}\alpha_i^2+\dfrac{1}{16}\sum_{i<j}\alpha_i\alpha_j-\dfrac{\rho}{2(1+2\rho)}
\sum_{\substack{i=1\\ mod 5}}^4\alpha_i\alpha_{i+1}.
\]

The infimum of $\dfrac{\Gamma_2(f)(v_0)}{\Gamma(f)(v_0)}$ is the least eigenvalue of the matrix
\[
\setlength\arraycolsep{8pt}
\begin{pmatrix}
\dfrac{\rho}{1+2\rho}+\dfrac{5}{4} & \dfrac{1}{4}-\dfrac{2\rho}{1+2\rho} & \dfrac{1}{4} & \dfrac{1}{4}-\dfrac{2\rho}{1+2\rho}\\[2em]
 \dfrac{1}{4}-\dfrac{2\rho}{1+2\rho} & \dfrac{\rho}{1+2\rho}+\dfrac{5}{4}  & \dfrac{1}{4}-\dfrac{2\rho}{1+2\rho} & \dfrac{1}{4}\\[2em]
  \dfrac{1}{4} & \dfrac{1}{4}-\dfrac{2\rho}{1+2\rho} & \dfrac{\rho}{1+2\rho}+\dfrac{5}{4}  & \dfrac{1}{4}-\dfrac{2\rho}{1+2\rho}\\[2em]
   \dfrac{1}{4}-\dfrac{2\rho}{1+2\rho} & \dfrac{1}{4} & \dfrac{1}{4}-\dfrac{2\rho}{1+2\rho}  & \dfrac{\rho}{1+2\rho}+\dfrac{5}{4}
\end{pmatrix}.\]

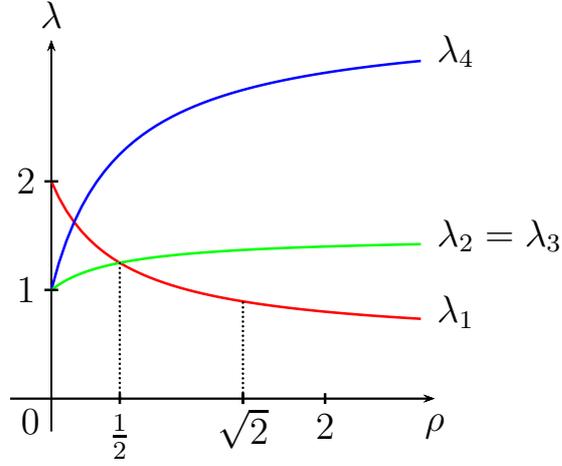
\begin{figure}[H]
\centering
\scalebox{1.2}{\begin{pspicture}*(-1,-1)(5.8,5.3)
\psset{xunit=1.5cm,yunit=1.2cm}
\SpecialCoor%
\psaxes[labels=none,ticks=none,linewidth=0.5pt]{->}(0,0)(-0.3,-0.3)(2.8,3.3)
\psPlot[linecolor=red]{0}{2.7}{(2+x)/(1+2*x)}
\psPlot[linecolor=blue]{0}{2.7}{(1+7*x)/(1+2*x)}
\psPlot[linecolor=green]{0}{2.7}{(1+3*x)/(1+2*x)}
\psline(1.4,-0.05)(1.4,0.05)
\psline[linestyle=dotted,dotsep=0.7pt](1.4,0)(1.4,0.89)
\psline[linestyle=dotted,dotsep=0.7pt](0.5,0)(0.5,1.25)
\uput[dl](0,0){0} \uput[d](2,0){2}  \psline(2,-0.05)(2,0.05)
\uput[d](0.5,0){$\frac{1}{2}$}  \psline(0.5,-0.05)(0.5,0.05)
\uput[d](1.4,0.05){$\sqrt{2}$}  \psline(1.4,-0.05)(1.4,0.05)
\uput[d](2.8,0){$\rho$} \uput[d](0,3.8){$\lambda$}
\uput[l](0,1){1}  \psline(-0.05,1)(0.05,1)
\uput[l](0,2){2}  \psline(-0.05,2)(0.05,2)
\uput[r](2.7,0.8){$\lambda_1$}
\uput[r](2.7,1.5){$\lambda_2=\lambda_3$}
\uput[r](2.7,3.2){$\lambda_4$}
\end{pspicture}}
\caption{The graph of the eigenvalues for $G_{4,\rho}$.}{\label{grafik2}}

\end{figure}

The eigenvalues are $\lambda_1=\dfrac{2+\rho}{1+2\rho}$, $\lambda_2=\lambda_3=\dfrac{1+3\rho}{1+2\rho}$ and $\lambda_4=\dfrac{1+7\rho}{1+2\rho}$, the least one being $\lambda_1$ on the interval $[1/2, \infty)$ (see Figure \ref{grafik2}) and $\lambda_2=\lambda_3$ on the interval $(0, 1/2]$.

\begin{figure}[H]
  \centering
  \includegraphics[scale=0.6]{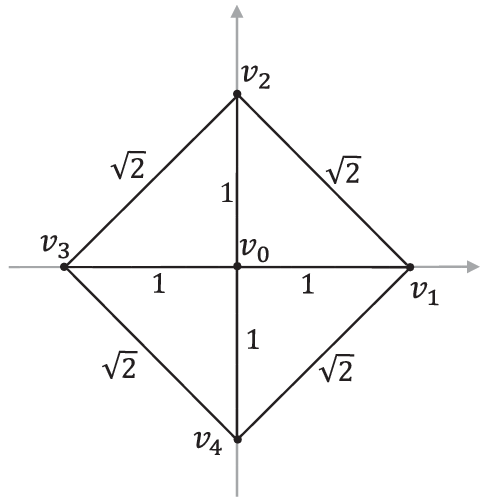}
   \caption{The umbrella graph $G_{4, \rho}$ with $\rho=\sqrt{2}$.} \label{sekil4}
\end{figure}

For $\rho=\sqrt{2}$, the umbrella graph $G_{4, \rho}$ is geodesically embeddable into the Euclidean plane (see Figure \ref{sekil4}) with Bakry-\'{E}mery curvature $K_G(v_0)=\dfrac{2+\sqrt{2}}{1+2\sqrt{2}}$.

\begin{figure}[H]
  \centering
  \includegraphics[scale=0.6]{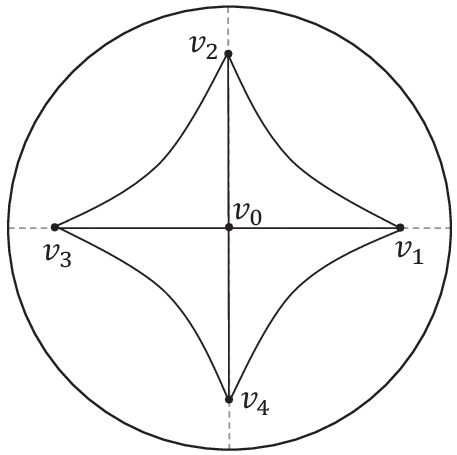}
   \caption{The umbrella graph $G_{4, \rho}$ with $\sqrt{2}<\rho<2$.} \label{sekil5}
\end{figure}

For $\sqrt{2}<\rho<2$, $G_{4, \rho}$ is geodesically embeddable into the hyperbolic plane with (smooth) curvature $-\dfrac{1}{K^2}$ (which is the solution of the equation $\cosh \dfrac{\rho}{K}=\cosh^2\dfrac{1}{K}$ ) (see Figure \ref{sekil5}) and has the discrete curvature $K_G(v_0)=\dfrac{2+\rho}{1+2\rho}$.

\begin{figure}[H]
  \centering
  \includegraphics[scale=0.6]{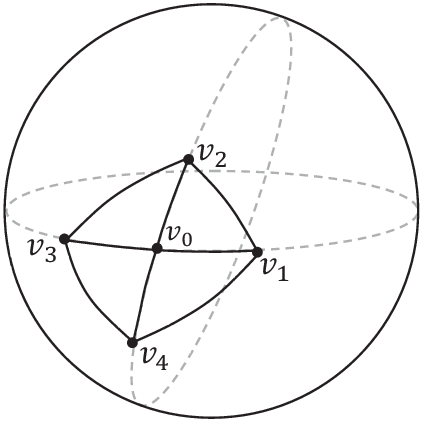}
   \caption{The umbrella graph $G_{4, \rho}$ with $0<\rho<\sqrt{2}$.} \label{sekil6}
\end{figure}

For $0<\rho<\sqrt{2}$, the graph $G_{4, \rho}$ is geodesically embeddable into the 2-sphere with radius $R$ (which is the solution of the equation $\cos \dfrac{\rho}{R}=\cos^2\dfrac{1}{R}$) and has the discrete curvature $\dfrac{2+\rho}{1+2\rho}$ for $1/2<\rho<\sqrt2$ and $\dfrac{1+3\rho}{1+2\rho}$ for $0<\rho<1/2$ (see Figure \ref{sekil6}).

We now give the results for $n=5$ and $n=6$, which are obtained with the help of the program in the Appendix.

\begin{figure}[H]
\centering
\scalebox{1.2}{\begin{pspicture}*(-1,-1)(5.8,4.8)
\psset{xunit=1.5cm,yunit=1.2cm}
\SpecialCoor%
\psaxes[labels=none,ticks=none,linewidth=0.5pt]{->}(0,0)(-0.3,-0.3)(2.8,2.8)
\psPlot[linecolor=red]{0}{2.7}{(2+x)/(1+2*x)}
\psPlot[linecolor=green]{0}{2.7}{(1+(4+5^(1/2))*x)/(1+2*x)}
\psPlot[linecolor=blue]{0}{2.7}{(1+(4-5^(1/2))*x)/(1+2*x)}
\psline(1.3,-0.05)(1.3,0.05)
\psline[linestyle=dotted,dotsep=0.7pt](1.3,0)(1.3,0.93)
\psline[linestyle=dotted,dotsep=0.7pt](1.175,0)(1.175,0.9)
\uput[dl](0,0){0} \uput[d](2,0){2} \psline(2,-0.05)(2,0.05)
\uput[d](1.05,0){$\rho_5^0$} \psline(1.175,-0.05)(1.175,0.05)
\uput[d](1.5,0){$\frac{1}{3-\sqrt{5}}$} \psline(1.3,-0.05)(1.3,0.05)
\uput[r](2.7,0.6){$\lambda_1$}
\uput[r](2.7,1.1){$\lambda_2=\lambda_3$}
\uput[r](2.7,2.8){$\lambda_4=\lambda_5$}
\uput[l](0,1){1}  \psline(-0.05,1)(0.05,1)
\uput[l](0,2){2}  \psline(-0.05,2)(0.05,2)
\uput[d](2.8,0){$\rho$} \uput[d](0,3.3){$\lambda$}
\end{pspicture}}
\caption{The graph of the eigenvalues for $G_{5,\rho}$. Here $\rho_5^0=\sqrt{3-\frac{1+\sqrt{5}}{2}}$.}\label{grafik3}
\end{figure}
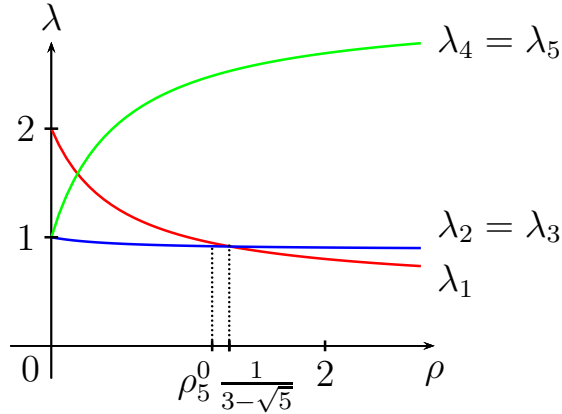

The auxiliary matrix $A_5$ for $n=5$ is
\[
\setlength\arraycolsep{8pt}
\begin{pmatrix}
 \dfrac{6+17\rho}{5(1+2\rho)} & \dfrac{1-8\rho}{5(1+2\rho)} & \dfrac{1}{5} & \dfrac{1}{5} & \dfrac{1-8\rho}{5(1+2\rho)} \\[2em]
  \dfrac{1-8\rho}{5(1+2\rho)} & \dfrac{6+17\rho}{5(1+2\rho)} & \dfrac{1-8\rho}{5(1+2\rho)} & \dfrac{1}{5} & \dfrac{1}{5} \\[2em]
 \dfrac{1}{5} & \dfrac{1-8\rho}{5(1+2\rho)} & \dfrac{6+17\rho}{5(1+2\rho)} & \dfrac{1-8\rho}{5(1+2\rho)} & \dfrac{1}{5}  \\[2em]
 \dfrac{1}{5} & \dfrac{1}{5}  & \dfrac{1-8\rho}{5(1+2\rho)} & \dfrac{6+17\rho}{5(1+2\rho)} & \dfrac{1-8\rho}{5(1+2\rho)} \\[2em]
 \dfrac{1-8\rho}{5(1+2\rho)} & \dfrac{1}{5}  & \dfrac{1}{5} & \dfrac{1-8\rho}{5(1+2\rho)} & \dfrac{6+17\rho}{5(1+2\rho)}
\end{pmatrix}\]
and its eigenvalues are $\lambda_1=\dfrac{2+\rho}{1+2\rho}$, $\lambda_2=\lambda_3=\dfrac{1+(4-\sqrt5)\rho}{1+2\rho}$ and $\lambda_4=\lambda_5=\dfrac{1+(4+\sqrt5)\rho}{1+2\rho}$, the least one being $\lambda_1$ on the interval $\left[\frac{1}{3-\sqrt5}, 2\right)$ and $\lambda_2$ on the interval $\left(0,\frac{1}{3-\sqrt5}\right]$ (see Figure \ref{grafik3}).

\begin{figure}[H]
\centering
\scalebox{1.2}{\begin{pspicture}*(-1,-1)(5.8,5)
\psset{xunit=1.5cm,yunit=1.2cm}
\SpecialCoor%
\psaxes[labels=none,ticks=none,linewidth=0.5pt]{->}(0,0)(-0.3,-0.3)(2.8,3.3)
\psPlot[linecolor=black]{0}{2.7}{(1+x)/(1+2*x)}
\psPlot[linecolor=red]{0}{2.7}{(2+x)/(1+2*x)}
\psPlot[linecolor=blue]{0}{2.7}{(1+5*x)/(1+2*x)}
\psPlot[linecolor=green]{0}{2.7}{(1+7*x)/(1+2*x)}
\psline(1,-0.05)(1,0.05)
\psline[linestyle=dotted,dotsep=0.7pt](1,0)(1,0.65)
\uput[dl](0,0){0} \uput[d](2,0){2} \psline(2,-0.05)(2,0.05)
\uput[d](1,0){$1$}  \psline(1,-0.05)(1,0.05)
\uput[l](0,1){1}  \psline(-0.05,1)(0.05,1)
\uput[l](0,2){2}  \psline(-0.05,2)(0.05,2)
\uput[d](2.8,0){$\rho$} \uput[d](0,3.7){$\lambda$}
\uput[r](2.7,0.9){$\lambda_1$} \uput[r](2.7,3.2){$\lambda_2$}
\uput[r](2.7,0.5){$\lambda_3=\lambda_4$}
\uput[r](2.7,2.28){$\lambda_5=\lambda_6$}
\end{pspicture}}
\caption{The graph of the eigenvalues for $G_{6,\rho}$. For the Euclidean plane embedding $\rho=1$.}\label{grafik4}

\end{figure}
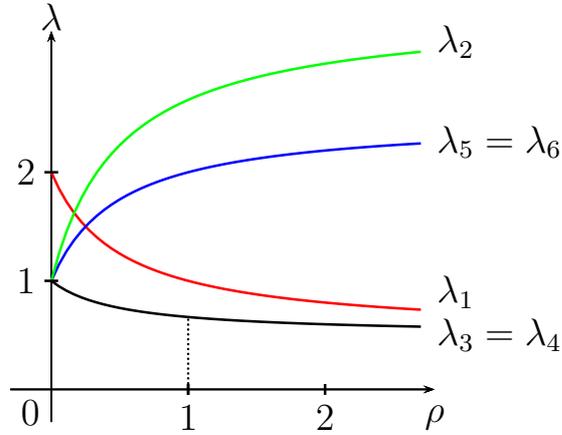

The auxiliary matrix $A_6$ for $n=6$ is
\[
\setlength\arraycolsep{8pt}
\begin{pmatrix}
 \dfrac{7+20\rho}{6(1+2\rho)} & \dfrac{1-10\rho}{6(1+2\rho)} &\dfrac{1}{6} &\dfrac{1}{6} & \dfrac{1}{6} & \dfrac{1-10\rho}{6(1+2\rho)} \\[2em]
 \dfrac{1-10\rho}{6(1+2\rho)} & \dfrac{7+20\rho}{6(1+2\rho)} & \dfrac{1-10\rho}{6(1+2\rho)} &\dfrac{1}{6} & \dfrac{1}{6}& \dfrac{1}{6} \\[2em]
 \dfrac{1}{6} &\dfrac{1-10\rho}{6(1+2\rho)} & \dfrac{7+20\rho}{6(1+2\rho)} & \dfrac{1-10\rho}{6(1+2\rho)} &\dfrac{1}{6} &\dfrac{1}{6} \\[2em]
 \dfrac{1}{6}&\dfrac{1}{6} & \dfrac{1-10\rho}{6(1+2\rho)} & \dfrac{7+20\rho}{6(1+2\rho)} & \dfrac{1-10\rho}{6(1+2\rho)} & \dfrac{1}{6}\\[2em]
\dfrac{1}{6} &\dfrac{1}{6} &\dfrac{1}{6} &\dfrac{1-10\rho}{6(1+2\rho)} & \dfrac{7+20\rho}{6(1+2\rho)} & \dfrac{1-10\rho}{6(1+2\rho)} \\[2em]
\dfrac{1-10\rho}{6(1+2\rho)} & \dfrac{1}{6} & \dfrac{1}{6} &\dfrac{1}{6} &\dfrac{1-10\rho}{6(1+2\rho)}  & \dfrac{7+20\rho}{6(1+2\rho)}
\end{pmatrix}.\]
and its eigenvalues are $\lambda_1=\dfrac{2+\rho}{1+2\rho}$, $\lambda_2=\dfrac{1+7\rho}{1+2\rho}$, $\lambda_3=\lambda_4=\dfrac{1+\rho}{1+2\rho}$, and $\lambda_5=\lambda_6=\dfrac{1+5\rho}{1+2\rho}$ with the least one being $\lambda_3=\lambda_4=\dfrac{1+\rho}{1+2\rho}$ (see Figure \ref{grafik4}).

Inspecting the graphs in Figures \ref{grafik1}, \ref{grafik2}, \ref{grafik3}, \ref{grafik4}, it is very remarkable to see that a ``spherical" n-umbrella graph (for $n=3, 4, 5, 6$) has always a higher discrete Bakry-\'{E}mery curvature than a ``hyperbolic" n-umbrella graph. For an arbitrary $n-$umbrella graph ($n\geq 7$), we cannot express the discrete curvature in closed form but with our calculator (given in the Appendix), we computed the values for $n=7, 8, 9, 10$ and $20$. The results are shown in Table \ref{table1}. For a fixed $n$ ($n=3,\dots,10, 20$), we first determine three values $\rho_n^0$,$\rho_n^+$ and $\rho_n^-$ of $\rho$ such that the $n-$umbrella graph $G_{n, \rho}$ can be geodesically embedded into the Euclidean plane for $\rho = \rho_n^0$, into the 2-sphere with radius 1 for $\rho = \rho_n^+$ and into the hyperbolic plane with curvature -1 for $\rho = \rho_n^-$. Then we compute the discrete Bakry-\'{E}mery curvature for these weighted graphs. This table shows that for each $n=3,\dots,10$ and $20$, the discrete curvatures for spherical cases are higher than the values of hyperbolic cases. We think that this surprising phenomenon is worth for further investigation.

\begin{table}[H]
  \centering
  \renewcommand{\arraystretch}{1.5}
  \begin{tabular}{ |c|c|c|c|c|c|c| }
 \hline
 $n$ & $\rho_n^+$ & $K_{G_{n,\rho_n^+}}(v_0)$ & $\rho_n^0$ & $K_{G_{n,\rho_n^0}}(v_0)$ & $\rho_n^-$ & $K_{G_{n,\rho_n^-}}(v_0)$ \\ \hline \hline
 3 & 1.6329 & 0.8516 & 1.7320 & 0.8360 & 1.7877 & 0.8278 \\ \hline
 4 & 1.2745 & 0.9226 & 1.4142 & 0.8918 & 1.5133 & 0.8725 \\ \hline
 5 & 1.0347 & 0.9204 & 1.1755 & 0.9171 & 1.2901 & 0.8918 \\ \hline
 6 & 0.8685 & 0.6826 & 1 & 0.6667 & 1.1163 & 0.6546 \\ \hline
 7 & 0.7474 & 0.5524 & 0.8677 & 0.5260 & 0.9800 & 0.5053 \\ \hline
 8 & 0.6557 & 0.4813 & 0.7653 & 0.4470 & 0.8716 & 0.4190 \\ \hline
 9 & 0.5835 & 0.4440 & 0.6840 & 0.4037 & 0.7836 & 0.3699 \\ \hline
 10 & 0.5261 & 0.4267 & 0.6180 & 0.3819 & 0.7112 & 0.3434 \\ \hline
 20 & 0.2640 & 0.5154 & 0.3128 & 0.4603 & 0.3656 & 0.4077 \\
 \hline
\end{tabular}
  \caption{$\rho_n^+$ (respectively $\rho_n^-$) is the value of $\rho$ for which the $n-$umbrella graph $G_{n,\rho}$ can be geodesically embedded into the sphere with radius 1 (respectively into the hyperbolic plane with curvature -1); $K_{G_{n,\rho}}(v_0)$ is the discrete Bakry-\'{E}mery curvature of $G_{n,\rho}$.}\label{table1}
\end{table}

\textbf{Declaration of Interest:} The authors report there are no competing interests to declare.

\newpage
\section{Appendix}

The following Maple code takes as input a finite, simple, weighted graph in form of a symmetric matrix $A$ with nonnegative entries. In ordering the vertices to build the matrix, take the vertex at which you want to compute the curvature in the first place. Then the program computes the discrete Bakry-\'{E}mery curvature at this vertex.

\begin{verbatim}
restart;
with(GraphTheory):
with(LinearAlgebra):
A := Matrix([ ]);    #Enter your graphs associated matrix here.
G := Graph(A);
N := RowDimension(A):
mu := MatrixVectorMultiply(A, Vector(N, 1)):
f := Vector(N, symbol = a):
f[1] := 0:
Komsu := Neighbors(G, 1):

DeltaMu := proc(fonk, x)
  local Komsux, v, j, toplam:
  toplam := 0:
  Komsux := Neighbors(G, x):
  for j to nops(Komsux) do
    v := Komsux[j]:
    toplam := toplam + A[x, v]*(fonk[v] - fonk[x]):
  end do:
  toplam := toplam/mu[x]:
  return toplam:
end proc:

GammaMu := proc(fonk1, fonk2, x)
  local Komsux, v, j, toplam:
  toplam := 0:
  Komsux := Neighbors(G, x):
  for j to nops(Komsux) do
    v := Komsux[j]:
    toplam := toplam + A[x, v]*(fonk1[v] - fonk1[x])*(fonk2[v] - fonk2[x]):
  end do:
  toplam := (1/(2*mu[x]))*toplam:
  return toplam:
end proc:

minimizeet := proc(x, w)
  local mindeger, Komsuxw, v, j, toplam1, toplam2:
  toplam1 := 0:
  toplam2 := 0:
  Komsuxw := convert(Neighbors(G, x), set)
                      intersect convert(Neighbors(G, w), set):
  Komsuxw := convert(Komsuxw, list):
  for j to nops(Komsuxw) do
    v := Komsuxw[j]:
    toplam1 := toplam1 + A[x, v]*A[v, w]*f[v]/mu[v]:
    toplam2 := toplam2 + A[x, v]*A[v, w]/mu[v]:
  end do:
  mindeger := simplify(2*toplam1/toplam2):
  return mindeger:
end proc:
DMF := Vector(N):
GMF := Vector(N):

for k to N do
    DMF[k] := DeltaMu(f, k):
end do:
for k to N do
    GMF[k] := GammaMu(f, f, k):
end do:
DMF:
GMF:
DeltaMu(GMF, 1):
GammaMu(f, DMF, 1):
GammaMu2X := 1/2*(DeltaMu(GMF, 1) - 2*GammaMu(f, DMF, 1)):
GammaMu2X := simplify(GammaMu2X):
GammaMu(f, f, 1):
derece := Degree(G, 1):
Pay := Matrix(derece):
Payda := Matrix(derece):
ikiKomsu := [];
for i to nops(Komsu) do
    ikiKomsu := {op(ikiKomsu), op(Neighbors(G, Komsu[i]))}:
end do:
convert(ikiKomsu, list):
T := convert(ikiKomsu, set):
birdisk := {1} union convert(Komsu, set):
ikiKomsu := T minus birdisk;
for i to nops(ikiKomsu) do
    w := ikiKomsu[i]:
    GammaMu2X := subs(f[w] = minimizeet(1, w), GammaMu2X):
end do:
simplify(expand(GammaMu2X)):
for i to derece do
    for j from i to derece do
      ii := Komsu[i]; jj := Komsu[j]:
      if ii = jj then polinom := GammaMu2X; polinom2 := GammaMu(f, f, 1):
        for k to derece do
          kk := Komsu[k]:
            if ii <> kk then polinom := subs(f[kk] = 0, polinom):
              polinom2 := subs(f[kk] = 0, polinom2):
            end if:
        end do:
        Pay[i, i] := simplify(subs(a[ii] = 1, polinom)):
        Payda[i, i] := simplify(subs(a[ii] = 1, polinom2)):
      end if:
      if ii <> jj then polinom := GammaMu2X:
        for k to derece do
          for m from k to derece do
            kk := Komsu[k]:
            mm := Komsu[m]:
            if ii <> kk and ii <> mm and jj <> kk and jj <> mm then
            polinom := subs(f[mm] = 0, polinom):
            end if:
          end do:
        end do:
        polinom := diff(polinom, a[ii]):
        polinom := diff(polinom, a[jj]):
        Pay[i, j] := simplify(polinom)/2:
        Pay[j, i] := simplify(polinom)/2:
      end if:
   end do:
end do:
Payda := MatrixInverse(Payda);
CM := MatrixMatrixMultiply(Payda, Pay);
EigenvaluesCM:= simplify(Eigenvalues(CM));
BEC := min(EigenvaluesCM);
\end{verbatim}

\end{document}